\RequirePackage{ifpdf}
\ifpdf 
\documentclass[pdftex]{sigma}
\else
\documentclass{sigma}
\fi

\DeclareMathAlphabet\oldmathcal{OMS}        {cmsy}{b}{n}
\SetMathAlphabet    \oldmathcal{normal}{OMS}{cmsy}{m}{n}
\DeclareMathAlphabet\oldmathbcal{OMS}       {cmsy}{b}{n}

\newcommand{\bb}{\mathbb}

\def\d{\partial}

\def\<{\langle}

\def\>{\rangle}

\def\BOne{{\mathchoice {\rm 1\mskip-4mu l} {\rm 1\mskip-4mu l}
                          {\rm 1\mskip-4.5mu l} {\rm 1\mskip-5mu l}}}

\def\fract#1#2{\raise4pt\hbox{$ #1 \atop #2 $}}
\def\decdnar#1{\phantom{\hbox{$\scriptstyle{#1}$}}
\left\downarrow\vbox{\vskip15pt\hbox{$\scriptstyle{#1}$}}\right.}
\def\decupar#1{\phantom{\hbox{$\scriptstyle{#1}$}}
\left\uparrow\vbox{\vskip15pt\hbox{$\scriptstyle{#1}$}}\right.}

\def\bbc{{\mathbb C}}

\def\bbp{{\mathbb P}}

\def\bbr{{\mathbb R}}

\def\bbz{{\mathbb Z}}

\def\gra{\alpha}

\def\gro{\omega}

\def\grt{\tau}

\def\grz{\zeta}
\def\grD{\Delta}

\def\grO{\Omega}

\def\bff{{\bf f}}

\def\bfz{{\bf z}}

\def\cala{{\mathcal A}}

\def\cald{{\mathcal D}}

\def\calf{{\mathcal F}}

\def\calh{{\mathcal H}}

\def\call{{\mathcal L}}

\def\cals{{\oldmathcal S}}

\def\la#1{\hbox to #1pc{\leftarrowfill}}
\def\ra#1{\hbox to #1pc{\rightarrowfill}}
\def\calz{{\oldmathcal Z}}

\def\X{\mathfrak{X}}

\def\ga{{\mathfrak a}}

\def\gc{{\mathfrak c}}

\def\ge{{\mathfrak e}}
\def\gf{{\mathfrak f}}
\def\gg{{\mathfrak g}}

\def\gi{{\mathfrak i}}

\def\gm{{\mathfrak m}}
\def\gn{{\mathfrak n}}
\def\go{{\mathfrak o}}

\def\gs{{\mathfrak s}}
\def\gt{{\mathfrak t}}

\def\gy{{\mathfrak y}}

\def\gB{{\mathfrak B}}
\def\gC{{\mathfrak C}}
\def\gD{{\mathfrak D}}

\def\gH{{\mathfrak H}}

\def\gS{{\mathfrak S}}

\def\hook{\mathbin{\hbox to 6pt{%
                 \vrule height0.4pt width5pt depth0pt
                 \kern-.4pt
                 \vrule height6pt width0.4pt depth0pt\hss}}}
\def\d{\partial}

\def\bJ{\bar{J}}

\def\bg{\bar{g}}
\def\bp{\bar{p}}

\def\tp{\tilde{p}}

\def\teta{\tilde{\eta}}

\def\con{{\gc\go\gn}}
\def\Con{{\gC\go\gn}}
\def\hX{\hat{X}}

\begin{document}

\allowdisplaybreaks

\renewcommand{\thefootnote}{$\star$}

\renewcommand{\PaperNumber}{058}

\FirstPageHeading

\ShortArticleName{Completely Integrable Contact Systems}

\AuthorNameForHeading{C.P.~Boyer}

\ArticleName{Completely Integrable Contact Hamiltonian Systems\\ and Toric Contact Structures on $\boldsymbol{S^2\times S^3}$\,\footnote{This paper is a
contribution to the Special Issue ``Symmetry, Separation, Super-integrability and Special Functions~(S$^4$)''. The
full collection is available at
\href{http://www.emis.de/journals/SIGMA/S4.html}{http://www.emis.de/journals/SIGMA/S4.html}}}

\Author{Charles P.~BOYER}

\Address{Department of Mathematics and Statistics, University of New Mexico,\\ Albuquerque, NM 87131, USA}
\Email{\href{mailto:cboyer@math.unm.edu}{cboyer@math.unm.edu}}
\URLaddress{\url{http://www.math.unm.edu/~cboyer/}}

\ArticleDates{Received January 28, 2011, in f\/inal form June 08, 2011;  Published online June 15, 2011}

\Abstract{I begin by giving a general discussion of completely integrable Hamiltonian systems in the setting of contact geometry. We then pass to the particular case of toric contact structures on the manifold $S^2\times S^3$. In particular we give a complete solution to the contact equivalence problem for a class of toric contact structures, $Y^{p,q}$, discovered by physicists
 by showing that $Y^{p,q}$ and $Y^{p',q'}$ are inequivalent as contact structures if and only if $p\neq p'$.}


\Classification{53D42; 53C25}

\begin{flushright}
\begin{minipage}{80mm}
\it Dedicated to Willard Miller Jr. on the occasion of his retirement
\end{minipage}
\end{flushright}


\renewcommand{\thefootnote}{\arabic{footnote}}
\setcounter{footnote}{0}

\section{Introduction}

This paper is based on a talk given at the S$^4$ conference at the University of Minnesota in honor of Willard Miller Jr. In turn that talk was based on my recent work in progress with J.~Pati~\cite{BoPa10} where we study the question of when certain toric contact structures on $S^3$-bundles over $S^2$ belong to equivalent contact structures. As in the talk, in this paper we concentrate on a particularly interesting special class of toric contact structures on $S^2\times S^3$ studied by physicists in~\cite{GMSW04a,MaSp05b,MaSp06}, and denoted by $Y^{p,q}$ where $p,q$ are relatively prime integers satisfying \mbox{$0<q<p$}. These structures have become of much interest in the study of the AdS/CFT conjecture \cite{GaMaSpWa04c,GaMaSpWa05} in M-theory since they admit Sasaki--Einstein metrics. The AdS/CFT correspondence relates string theory on the product of anti-deSitter space with a compact Einstein space to quantum f\/ield theory on the conformal boundary, thus giving a kind of holographic principle. As Sasaki--Einstein metrics admit Killing spinors \cite{FrKa89,FrKat2}, the string theories or M-theory are supersymmetric. The relation to contact structures is that Sasakian metrics are a special class of contact metric structures, and roughly Sasakian geometry is to contact geometry what K\"ahlerian geometry is to symplectic geometry. We refer to the recent book~\cite{BG05} for a thorough treatment of Sasakian geometry.

The connection between completely integrable Hamiltonian systems and toric geometry in the symplectic setting is best described by the famous Arnold--Liouville theorem\footnote{A very nice treatment is given by Audin~\cite{Aud08}.} which in its modern formulation (due to Arnold \cite{Arn78}) roughly states the following: let $(M^{2n},\gro)$ be a~symplectic manifold of dimension $2n$ with a Hamiltonian $h$, and assume that there are $n$ f\/irst integrals $\bff=(h=f_1,\dots,f_n)$ in involution that are functionally independent on a dense open subset of~$M$. Such a structure is called a {\it completely integrable Hamiltonian system}. Let $a$ be a regular value of the {\it moment map} $\bff:M\ra{1.6} \bbr^n$, and assume that the f\/iber $\bff_a=\bff^{-1}(a)$ is compact and connected, then $\bff_a$ is a torus $T^n$, and moreover, there is a neighborhood of $\bff_a$ that is dif\/feomorphic to $T^n\times D^n$ where $D^n$ is an $n$-dimensional disk, and the f\/low of $h$ is linear in the standard coordinates on~$T^n$ and independent of the coordinates of $D^n$. The coordinates of $T^n$ are called {\it angle coordinates} and those of $D^n$ {\it action coordinates}. Thus, locally such a manifold looks like a toric symplectic manifold, that is, a symplectic manifold with a locally free local torus action. However, there is an obstruction to having a global torus action~\cite{Dui80,BoMo89}, namely the monodromy of a certain period lattice. The case where one does have a global Hamiltonian $T^n$-action on a compact symplectic manifold $(M^{2n},\gro)$ is both beautiful and well-understood. First, there is the Atiyah--Guillemin--Sternberg theorem \cite{At82,GuSt82} which says that the image of the moment map is a convex polytope in $\bbr^n$, and then a theorem of Delzant \cite{Del88} which states that the polytope characterizes the toric symplectic structure up to equivariant Hamiltonian symplectomorphism.

Turning to the contact case, the development has been more recent. In fact, developing a theory of completely integrable systems in contact geometry was listed as problem \#1995-12 in~\cite{Arn04}.
Arnold seemed to have been unaware of the seminal work of Banyaga and Molino~\cite{BM93,BM96,Ban99} who develop the case of a local action of an $(n+1)$-dimensional torus on an oriented compact contact manifold giving the contact version of the Arnold--Liouville theorem under some additional assumptions. But even a bit earlier the foliation approach to contact complete integrability was given~\cite{Lib91,Pan90}. Much more recently a description in terms of a f\/lag of foliations was given in~\cite{KhTa09}. The approach presented here is more along the classical lines of using f\/irst integrals of commuting functions. As we shall see there are some subtle dif\/ferences with the symplectic case which manifest themselves dif\/ferently depending on the presentation. As mentioned above our main focus will be on completely integrable contact systems on a $(2n+1)$-dimensional compact contact manifold that arise from the global action of an $(n+1)$-dimensional torus.

As in the symplectic case the monodromy of an appropriate period lattice is the obstruction to having a global $T^{n+1}$-action. In \cite{BG00b} the subclass of contact manifolds with a $T^{n+1}$-action whose Reeb vector f\/ield lies in the Lie algebra of the torus was studied. It was shown that all such toric contact manifolds (of {\it Reeb type}) are determined by a certain polytope lying in a hyperplane (the {\it characteristic hyperplane}) in the dual of the Lie algebra of the torus, and they can all be obtained from contact reduction of an odd dimensional sphere with its standard contact structure. Furthermore, all toric contact structures of Reeb type admit a compatible Sasakian metric. A complete classif\/ication of all compact toric contact manifolds up to $T^{n+1}$-equivariance was then given by Lerman~\cite{Ler02a}. We are interested in the contact equivalence problem in the toric setting. We can ask the following question. Given any two inequivalent toric contact Hamiltonian structures on a smooth manifold $M$, when are they equivalent as contact manifolds? Although there are several new results in this paper, its main purpose is to give a proof of the following theorem which is a particular case of the more general results to appear in~\cite{BoPa10}.

\begin{theorem}\label{SEequiv}
Let $p,q$ be relatively prime integers satisfying $0<q<p$.
The toric contact structures $Y^{p,q}$ and $Y^{p',q'}$ on $S^2\times S^3$ belong to equivalent contact structures if and only if $p'=p$, and for each fixed integer $p>1$ there are exactly $\phi(p)$ toric contact structures $Y^{p,q}$ on $S^2\times S^3$ that are equivalent as contact structures. Moreover, the contactomorphism group of $Y^{p,q}$ has at least $\phi(p)$ conjugacy classes of maximal tori of dimension three.
\end{theorem}

Here $\phi(p)$ denotes the Euler phi function, that is the number of positive integers that are less than $p$ and relatively prime to $p$.

\section{A brief review of contact geometry}

In this section we give a very brief review of contact geometry referring to the books \cite{LiMa87,KrVi99,BG05,Bla10} for details.

\subsection{Contact manifolds}\label{section1.1}

Recall that a {\it contact structure} on a connected oriented manifold $M$ is an equivalence class of 1-forms $\eta$ satisfying $\eta\wedge (d\eta)^n\neq 0$ everywhere on $M$ where two 1-forms $\eta$, $\eta'$ are equivalent if there exists a nowhere vanishing function $f$ such that $\eta'=f\eta$. We shall also assume that our contact structure has an orientation, or equivalently, the function $f$ is everywhere positive. More conveniently the contact structure can be thought of as the oriented $2n$-plane bundle def\/ined by $\cald=\ker\eta$. A manifold $M$ with a contact structure $\cald$ is called a {\it contact manifold} which is necessarily odd dimensional, and is denoted by $(M,\cald)$. Choosing a contact form $\eta$ gives $\cald$ the structure of a symplectic vector bundle with 2-form $d\eta$. Choosing another contact form $\eta'=f\eta$ we see that
\begin{gather}\label{confsymp}
d\eta'|_{\cald\times\cald}=fd\eta|_{\cald\times\cald},
\end{gather}
so $\cald$ has a natural conformal symplectic structure.

For every choice of contact 1-form $\eta$ there exists a unique vector f\/ield $R_\eta$, called the {\it Reeb vector field}, that satisf\/ies $\eta(R_\eta)=1$ and $R_\eta\hook d\eta=0$. The dynamics of the Reeb f\/ield $R_\eta$ can change drastically as we change $\eta$. The one dimensional foliation $\calf_{R_\eta}$ on $M$ generated by~$R_\eta$ is often called the {\it characteristic foliation}. We say that the foliation $\calf_{R_\eta}$ is {\it quasi-regular} if there is a positive integer~$k$ such that each point has a
foliated coordinate chart $(U,x)$ such that each leaf of
$\calf_{R_\eta}$ passes through $U$ at most $k$ times. If $k=1$ then
the foliation is called {\it regular}. We also say that the corresponding contact 1-form $\eta$ is {\it quasi-regular $($regular$)$}, and more generally that a contact structure $\cald$ is {\it quasi-regular $($regular$)$} if it has a quasi-regular (regular) contact 1-form. A contact 1-form (or characteristic foliation) that is not quasi-regular is called {\it irregular}.
When $M$ is compact a regular contact form $\eta$ is a connection 1-form in a principle $S^1$ bundle $\pi:M\ra{1.5} \calz$ over a symplectic base manifold $\calz$ whose symplectic form $\gro$ satisf\/ies $\pi^*\gro=d\eta$. In the quasi-regular case $\pi:M\ra{1.5} \calz$ is an $S^1$ orbibundle over the symplectic orbifold $\calz$. The former is known as the {\it Boothby--Wang construction} \cite{BoWa} and the latter the {\it orbifold Boothby--Wang construction}~\cite{BG00a}. $S^1$ orbibundles play an important role in the proof of Theorem~\ref{SEequiv}.

\subsection{Compatible metrics and Sasakian structures}

Let $(M,\cald)$ be a contact manifold and f\/ix a contact form $\eta$. Choose an almost complex structure~$J$ in the symplectic vector bundle $(\cald,d\eta)$ and extend it to a section $\Phi$ of the endomorphism bundle of $TM$ by demanding that it annihilates the Reeb vector f\/ield, that is, $\Phi R_\eta=0$. We say that the almost complex structure~$J$ is {\it compatible} with $\cald$ if for any sections $X$, $Y$ of $\cald$ we have
\[
d\eta(JX,JY)=d\eta(X,Y), \qquad d\eta(JX,Y)>0.
\]
Note that $g_\cald(X,Y)=d\eta(JX,Y)$ def\/ines an Hermitian metric on the vector bundle $\cald$. Moreover, we can extend this to a Riemannian metric on $M$ by def\/ining
\begin{gather*}
g=d\eta\circ(\Phi\otimes \BOne)+\eta\otimes \eta.
\end{gather*}
Note that the contact metric $g$ satisf\/ies the compatibility condition
\[
g(\Phi X,\Phi Y)=g(X,Y) -\eta(X)\eta(Y),
\]
where $X$, $Y$ are vector f\/ields on $M$.
Then the quadruple $\cals=(R_\eta,\eta,\Phi,g)$ is called a {\it contact metric structure} on $M$. Note also that the pair $(\cald,J)$ def\/ines a strictly pseudoconvex almost CR structure on $M$. The contact metric structure $(R_\eta,\eta,\Phi,g)$ is said to be {\it K-contact} if the Reeb vector f\/ield $R_\eta$ is a Killing vector f\/ield for the metric $g$, that is, if $\pounds_{R_\eta}g=0$. This is equivalent to the condition $\pounds_{R_\eta}\Phi=0$. If in addition the almost CR structure $(\cald,J)$ is integrable, that is a CR structure, then $(R_\eta,\eta,\Phi,g)$ is a~{\it Sasakian structure}. For a detailed treatment, including many examples, of Sasakian structures we refer to~\cite{BG05}.

If $M$ is compact and $\cals=(R_\eta,\eta,\Phi,g)$ is a Sasakian structure (actually K-contact is enough) on $M$, then if necessary by perturbing the Reeb vector f\/ield we can take $R_\eta$ to generate an $S^1$-action which leaves invariant the Sasakian structure. So the Sasakian automorphism group ${\rm Aut}(\cals)$ which is a compact Lie group has dimension at least one. If its maximal torus $T^k$ has dimension $k$ greater than one, then there is a cone $\gt^+_k$ of Reeb vector f\/ields, the {\it Sasaki cone}, lying in the Lie algebra $\gt_k$ of $T^k$ such that $\eta(\xi)>0$ everywhere for all $\xi\in \gt_k^+$. Note that the vector f\/ield $\xi$ is the Reeb vector f\/ield for the contact form $\eta'= {\frac{\eta}{\eta(\xi)}}$, and the induced contact metric structure $\cals'=(\xi,\eta',\Phi',g')$ is Sasakian. The conical nature of $\gt_k$ is exhibited by the {\it transverse homothety} (cf.~\cite{BG05}) which takes a Sasakian structure $\cals=(\xi,\eta,\Phi,g)$ to the Sasakian structure
\[
\cals_a=\big(a^{-1}\xi,a\eta,\Phi,ag+\big(a^2-a\big)\eta\otimes \eta\big)
\]
for any $a\in \bbr^+$.

\subsection{The symplectization}

Contact geometry can be understood in terms of symplectic geometry through its symplectization. Given a contact structure $\cald$ on~$M$ we recall the symplectic cone $C(M)=M\times \bbr^+$ with its natural symplectic structure $\grO=d(r^2\eta)$ where $r$ is a coordinate on $\bbr^+$. Note that~$\grO$ only depends on the contact structure $\cald$ and not on the choice of contact form~$\eta$. For if $\eta'=e^{2f}\eta$ is another choice of contact form, we can change coordinates $r'=e^{-f}r$ to give $d(r'^2\eta')=d(r^2\eta)=\grO$. The symplectic cone $(C(M),\grO)$ is called the {\it symplectization} or the symplectif\/ication of $(M,\cald)$. Recall the Liouville vector f\/ield $\Psi=r\frac{\d}{\d r}$ on the cone $C(M)$ and notice that it is invariant under the above change of coordinates, i.e., $\Psi=r\frac{\d}{\d r}=r'\frac{\d}{\d r'}.$
We have chosen the dependence of $\grO$ on the radial coordinate to be homogeneous of degree $2$ with respect to $\Psi$, since we want compatibility with cone metrics and these are homogeneous of degree $2$. In fact, a contact metric structure $(R_\eta,\eta,\Phi,g)$ on $M$ gives rise to an almost K\"ahler structure $(\grO,\bg=dr^2+r^2g)$ on $C(M)$ which is K\"ahler if and only if $(R_\eta,\eta,\Phi,g)$ is Sasakian.

An alternative approach to the symplectization is to consider the cotangent bundle $T^*M$ with its canonical (tautological) 1-form def\/ined as follows. It is the unique 1-form $\theta$ on $T^*M$ such that for every 1-form $\gra:M\ra{1.8} T^*M$ we have $\gra^*\theta=\gra$.  In local coordinates $(x^i,p_i)$ on $T^*M$ the canonical 1-form is given by $\theta=\sum_ip_i dx^i$. This gives $T^*M$ a canonical symplectic structure def\/ined by $d\theta$. Let $\cald^o$ be the annihilator of $\cald$ in $T^*M$ which is a real line bundle on $M$, and a~choice of contact 1-form $\eta$ trivializes $\cald^o\approx M\times \bbr$. Then $\cald^o\setminus \{0\}$ splits as $\cald^o\setminus \{0\}\approx \cald^o_+\cup \cald^o_-$, where the sections of $\cald^o_+$ are of the form $f\eta$ with $f>0$ everywhere on $M$. Thus, we have the identif\/ication $C(M)= M\times \bbr^+\approx \cald^o_+$ which is also identif\/ied with the principal $\bbr^+$ bundle associated to the line bundle $\cald^o$. From a more intrinsic viewpoint the symplectization is the total space of the principal $\bbr^+$-bundle $\cald^o_+$. A choice of oriented contact form $\eta$ gives a global section of $\cald^o_+$, and hence a trivialization of $\cald^o_+$. Now $\teta=r^2\eta$ is a 1-form on $C(M)$, so $\teta^*\theta=\teta$. Thus, the symplectic form $\grO$ on $C(M)$ satisf\/ies $\grO=d\teta=\teta^*d\theta$.

\subsection{The group of contactomorphisms}\label{congroupsect}

We are interested in the subgroup $\gC\go\gn(M,\cald)$ of the group $\gD\gi\gf\gf(M)$ of all dif\/feomorphisms of $M$ that leave the contact structure $\cald$ invariant. Explicitly, this {\it contactomorphism group} is def\/ined by
\begin{gather*}
\gC\go\gn(M,\cald)=\{\phi\in \gD\gi\gf\gf(M)~|~\phi_*\cald\subset \cald\}.
\end{gather*}
We are actually interested in the subgroup $\gC\go\gn(M,\cald)^+$ of contactomorphisms that preserve the orientation of $\cald$. Alternatively, if we choose a contact form $\eta$ representing $\cald$ these groups can be characterized as
\begin{gather*}
\gC\go\gn(M,\cald) =\{\phi\in \gD\gi\gf\gf(M)~|~\phi^*\eta=f\eta,~f(x)\neq 0~\text{for all}~x\in M\},\\
\gC\go\gn(M,\cald)^+ =\{\phi\in \gD\gi\gf\gf(M)~|~\phi^*\eta=f\eta,~f(x)>0~\text{for all}~x\in M\}.
\end{gather*}

We are also mainly concerned with the case that the manifold $M$ is compact. In this case $\gD\gi\gf\gf(M)$ and $\gC\go\gn(M,\cald)$ can be given the compact-open $C^\infty$ topology\footnote{Generally, in the non-compact case, this topology does not control the behavior at inf\/inity, and a much larger topology should be used.} in which case $\gC\go\gn(M,\cald)$ becomes a regular Fr\'echet Lie group \cite{Omo97,Ban97} locally modelled on the Fr\'echet vector space $\gc\go\gn(M,\cald)$ of inf\/initesimal contact transformations, that is the Lie algebra of $\gC\go\gn(M,\cald)$ def\/ined by
\begin{gather*}
\gc\go\gn(M,\cald)=\{X\in \X(M)~|~\text{if}~Y~\text{is a $C^\infty$ section of}~\cald,~\text{so is}~[X,Y]\},
\end{gather*}
where $\X(M)$ denotes the vector space of all $C^\infty$ vector f\/ields on $M$. It is easy to see that this is equivalent to the condition
\begin{gather}\label{conliealg}
\pounds_X\eta=a_X\eta
\end{gather}
for any contact form $\eta$ representing $\cald$ and some $a_X\in C^\infty(M)$. We are also interested in the subgroup $\gC\go\gn(M,\eta)$ consisting of all $\phi\in \gC\go\gn(M,\cald)$ such that $\phi^*\eta=\eta$. Its Lie algebra is
\begin{gather*}
\gc\go\gn(M,\eta)=\{X\in \gc\go\gn(M,\cald)~|~\pounds_X\eta=0\}.
\end{gather*}

Similarly, on a symplectic manifold $(N,\gro)$ we have the group $\gS\gy\gm(N,\gro)$ of symplectomorphisms def\/ined by
\begin{gather*}
\gS\gy\gm(N,\gro)=\{\phi\in \gD\gi\gf\gf(N)~|~\phi^*\gro=\gro\}.
\end{gather*}
When $N$ is compact this group is also a regular Fr\'echet Lie group locally modelled on its Lie algebra
\begin{gather*}
\gs\gy\gm(N,\gro)=\{X\in \X(M)~|~\pounds_X\gro=0\}.
\end{gather*}
However, we are interested in the symplectic cone $(C(M),\grO)$ which is non-compact. Fortunately, the subgroup of $\gS\gy\gm(C(M),\grO)$ that is important for our purposes behaves as if $C(M)$ were compact. Let $\gD=\gD(C(M))$ denote the 1-parameter group of dilatations generated by the Liouville vector f\/ield $\Psi$, and let $\gS\gy\gm(C(M),\grO)^\gD$ denote the subgroup  consisting of all ele\-ments of $\gS\gy\gm(C(M),\grO)$ that commute with $\gD$. Then one easily sees \cite{LiMa87} that there is an isomorphism of groups $\gS\gy\gm(C(M),\grO)^\gD\approx \gC\go\gn(M,\cald)$. On the inf\/initesimal level we also have $\gc\go\gn(M,\cald)\approx \gs\gy\gm(C(M),\grO)^\gD$ where $\gs\gy\gm(C(M),\grO)^\gD$ denotes the Lie subalgebra of all elements of $\gs\gy\gm(C(M),\grO)$ that commute with $\Psi$. This isomorphism is given explicitly by
\begin{gather*}
X\mapsto X-\frac{a_X}{2}\Psi.
\end{gather*}
Since $a_X$ is def\/ined by equation~(\ref{conliealg}) this isomorphism depends on the choice of contact form $\eta$.

\subsection{Legendrian and Lagrangian submanifolds}

Recall that a subspace $E$ of a symplectic vector space $(V,\gro)$ is {\it isotropic} or {\it $($co-isotropic$)$} if $E\subset E^\perp$ or $(E^\perp\subset E)$, respectively, where $E^\perp$ denotes the symplectic orthogonal to $E$. A~{\it Lagrangian subspace} is a maximal isotropic subspace or equivalently one which is both isotropic and co-isotropic, i.e.~$E=E^\perp$. A submanifold $f:P\ra{1.8} N$ of a symplectic manifold $(N^{2n},\gro)$ whose tangent space at each point $p\in P$ is a Lagrangian subspace of $(f^*TN)_p$ with respect to $f^*\gro$ is called a {\it Lagrangian submanifold}, and has dimension $n$. Locally all symplectic manifolds look the same, and so do all Lagrangian submanifolds. In a local Darboux coodinate chart $(p_i,q^i)$ we have $\gro=\sum_idp_i\wedge dq^i$, and the Lagrangian submanifolds are the leaves of a foliation, called the {\it Lagrangian foliation}, generated by the vector f\/ields $\{\d_{p_i}\}_{i=1}^n$. These vector f\/ields form an $n$-dimensional Abelian subalgebra of $\gs\gy\gm(N,\gro)$.

Now consider the case of a contact manifold $(M^{2n+1},\cald)$ with its natural conformal symplectic structure described by equation~(\ref{confsymp}). Then the isotropic and co-isotropic subspaces of $\cald_p$ at a point $p\in M$ are independent of the choice of $\eta$, and the maximal isotropic (Lagrangian) subspaces have dimension $n$. An integral submanifold of $\cald$ whose tangent spaces are Lagrangian subspaces of $\cald$ is called a {\it Legendrian submanifold}. As in the symplectic case locally all contact manifolds are the same, and any contact 1-form $\eta$ can be written in a Darboux coordinate chart $(z,p_i,q^i)$ as $\eta=dz-\sum_ip_idq^i$. Again the vector f\/ields $\{\d_{p_i}\}_{i=1}^n$ give a foliation, called a~{\it Legendrian foliation} whose leaves are Legendrian submanifolds; however, these vector f\/ields are {\it not} inf\/initesimal contact transformation, since
\[
\eta\wedge\pounds_{\d_{p_i}}\eta =\eta\wedge \bigl(\d_{p_i}\hook d\eta+d(\eta(\d_{p_i})\bigr) =-\eta\wedge dq^i\neq 0.
\]
Nevertheless, one easily f\/inds inf\/initesimal contact transformations whose projections onto~$\cald$ are~$\d_{p_i}$, namely the vector f\/ields $\d_{p_i}+q^i\d_z$. These generate an $n$-dimensional Abelian subalgebra of $\gc\go\gn(M,\eta)$. Note that the Reeb vector f\/ield of $\eta$ is $\d_z$, and that the vector f\/ields $\{\d_{p_i}+q^i\d_z,\d_z\}_{i=1}^n$ span an $(n+1)$-dimensional Abelian Lie algebra of $\gc\go\gn(M,\eta)$, and describes the `co-Legendrian foliation' of \cite{KhTa09}. Actually, the local geometry of a contact structure is described by the vector f\/ields $\{\d_{p_i}+q^i\d_z,\d_{q^i},\d_z\}_{i=1}^n$ which span the Lie algebra of the Heisenberg group $\calh^{2n+1}\subset \Con(\bbr^{2n+1},\eta)$.

The main interest in Lagrangian and Legendrian submanifolds is with their global behavior. Moreover, generally they do not form foliations, but {\it singular foliations}, and the nature of the singularities are often related to the topology of the underlying manifold.

Let us relate Legendrian submanifolds of a contact manifold $(M,\cald)$ to Lagrangian sub\-mani\-folds of the symplectization $(C(M),\grO)$. Choosing a contact form we can write
\begin{gather*}
\grO=d\big(r^2\eta\big)=2rdr\wedge \eta +r^2d\eta.
\end{gather*}
So we see that a Legendrian submanifold $\call$ of $(M,\cald)$ lifts to an isotropic submanifold $\tilde{\call}$ of $(C(M),\grO)$ which by equation~(\ref{confsymp}) is independent of the choice of contact form $\eta$. Since Lagrangian submanifolds $L$ of $(C(M),\grO)$ have dimension $n+1$, the lift $\tilde{\call}$ is a codimension one submanifold of some $L$.

\subsection{Invariants of contact structures}\label{conhom}

It is well known that as in symplectic geometry there are no local invariants in contact geometry. Indeed, if $(\cald_t,\eta_t)$ denotes a 1-parameter family of contact structures on $M$ with $t\in [0,1]$, then Gray's theorem says that there exists a dif\/feomorphism $\varphi_t:M\ra{1.6} M$ such that $\varphi_t^*\eta_t=f_{\varphi_t}\eta_0$ for each $t\in [0,1]$. The simplest invariant is the f\/irst Chern class $c_1(\cald)$ of the symplectic vector bundle $(\cald,d\eta)$. Note two remarks: f\/irst since the set of isomorphism classes of symplectic vector bundles coincides with the set of isomorphism classes of complex vector bundles, the Chern classes of $\cald$ are well def\/ined; second, $c_1(\cald)$ is independent of the choice of contact form~$\eta$ since if $\eta'=f\eta$ for some nowhere vanishing smooth function on $M$, then equation (\ref{confsymp}) holds. So if $c_1(\cald')\neq c_1(\cald)$, then the contact structures $\cald$ and $\cald'$ are inequivalent.

However, all the contact structures in Theorem~\ref{SEequiv} have $c_1(\cald)=0$. Thus, it is important to distinguish contact structures with the same f\/irst Chern class. Fortunately, this can be done with contact homology which is a piece of the larger symplectic f\/ield theory of Eliashberg, Givental, and Hofer \cite{ElGiHo00}. We do not go into details here as it would take us too far af\/ield, but only sketch the idea and refer to the literature \cite{Bou02t,Bou03,ElGiHo00} for details. The idea is to construct a Floer-type homology theory on the free loop space of a contact manifold $M$. Fix a contact form $\eta$ and consider the action functional $\mathcal{A}: C^{\infty}(S^1; M) \to \bb{R},$
def\/ined by
\begin{gather*}
\cala(\gamma) = \int_{\gamma} \eta.
\end{gather*}
The critical points of $\mathcal{A}$ are closed orbits of the Reeb vector f\/ield of $\eta,$ and the gradient trajectories, considered
as living in the symplectization $C(M)$ of $M$ are pseudoholomorphic curves which are cylindrically asymptotic over closed Reeb orbits. The idea then is to construct a chain complex $C_{*}$ generated by closed Reeb orbits of a suitably generic Reeb vector f\/ield. The homology of this complex is called {\it contact homology}, and its grading is determined by the Conley--Zehnder index (or Robbin--Salamon index) which roughly speaking measures the twisting of nearby Reeb orbits about a closed Reeb orbit. The dif\/ferential in this complex is given by a suitable count of pseudoholomorphic curves in the symplectization $C(M)$ of $M$. Assuming a certain transversality condition holds\footnote{A full treatment of the transversality problem awaits the completion of polyfold theory by Hofer, Wysocki, and Zehnder. See for example \cite{Hof08p} and references therein.}, the contact homology ring $HC_*$ is an invariant of the contact structure, and thus can be used to distinguish contact structures with the same f\/irst Chern class.

Another invariant of a contact manifold $(M^{2n+1}\!,\cald)$ is its contactomorphism group $\Con(M,\cald)$; however, it is too big to be of much use. On the other hand the number $\gn(\cald,k)$ of conjugacy classes of maximal tori of dimension $k\leq n+1$ in $\Con(M,\cald)$ is also an invariant, so it too can be used to distinguish contact structures. The problem here is that unlike the symplectomorphism group it is dif\/f\/icult to get a precise answer for $\gn(\cald,k)$. In our proof of Theorem~\ref{SEequiv} given in Section~\ref{proofthm} we can only obtain a lower bound for the number of conjugacy classes of $3$-tori, namely $\gn(\cald_p,3)\geq \phi(p)$.

\section{Toric contact structures as completely integrable Hamiltonian systems}

While completely integrable Hamiltonian systems in symplectic geometry have a long and distinguished history, as mentioned previously the contact version of complete integrability has been considered only fairly recently \cite{BM93,BM96,BG00b,Ler02a,KhTa09}. In symplectic geometry complete integrability can be def\/ined in terms of certain possibly singular foliations called Lagrangian foliations, and the vectors tangent to the leaves of this foliation provide an Abelian Lie algebra of inf\/initesimal symplectic transformations of one half the dimension of the manifold at least locally. The si\-tuation in contact geometry is somewhat more subtle. While the contact bundle $\cald$ does have a similar foliation, the Legendre foliation, its generic leaves have dimension one less than complete integrability requires. Moreover, non-trivial sections of $\cald$ are never inf\/initesimal contact transformations, so one must extend these sections f\/irst. Then one could just add in the Reeb vector f\/ield to obtain an Abelian Lie algebra of the correct dimension. However, this falls short of capturing all cases.

Recall that in symplectic geometry a Hamiltonian for the symplectic structure $\gro$ is a smooth function $H$ such that $X\hook \gro=-dH$ where $X$ is an inf\/initesimal symplectomorphism, that is, $\pounds_X\gro=0$. However, such an $H$ exists only if the de Rham cohomology class of $X\hook \gro$ vanishes. When it does exist the corresponding vector f\/ield $X_H$ is called a  {\it Hamiltonian vector field} and the triple $(N,\gro,X_H)$ is called a {\it Hamiltonian system}. Unlike the symplectic case, contact structures are automatically Hamiltonian; however, one needs to choose a certain isomorphism as described below.

\subsection{Contact Hamiltonian systems}
It is a well known result of Libermann (cf.~\cite{LiMa87}) that a choice of contact 1-form $\eta$ gives an isomorphism between the Lie algebra of inf\/initesimal contact transformations $\gc\go\gn(M,\cald)$ and the Lie algebra of smooth functions $C^\infty(M)$ by sending $X\in \gc\go\gn(M,\cald)$ to $\eta(X)\in C^\infty(M)$. The Lie algebra structure on $C^\infty(M)$ is given by the Jacobi bracket $\{\eta(X),\eta(Y)\}_\eta=\eta([X,Y])$. We then call the function $\eta(X)$ the {\it contact Hamiltonian} associated to the contact vector f\/ield $X$.
So any smooth function $f$ on $M$ can be a contact Hamiltonian, but it entails a choice of contact form. Moreover, as indicated the Jacobi bracket itself depends on the choice of contact form. Let $\eta'=f\eta$ be another contact form compatible with the co-orientation, so $f>0$ everywhere, and let $g,h\in C^\infty(M)$. Then the corresponding Jacobi brackets are related by
\begin{gather*}
\{g,h\}_{\eta'}=f\left\{\frac{g}{f},\frac{h}{f}\right\}_\eta.
\end{gather*}

Note that unlike the Poisson bracket, the Jacobi bracket does not satisfy the Leibniz rule, and $\{g,1\}_\eta=0$ if and only if $[X_g,R_\eta]=0$ where $R_\eta$ is the Reeb vector f\/ield of $\eta$, and $X_g$ is the Hamiltonian vector f\/ield corresponding to $g$. Furthermore, it is well known \cite{LiMa87,BG05} that the centralizer of $R_\eta$ in $\gc\go\gn(M,\cald)$ is the Lie subalgebra $\gc\go\gn(M,\eta)$. If we f\/ix a contact form $\eta$ and consider a Hamiltonian $h_X=\eta(X)$, we can contract equation (\ref{conliealg}) with the Reeb vector f\/ield~$R_\eta$ of~$\eta$ and use the well known Cartan equation $\pounds_X\eta =X\hook d\eta +dh_X$ to give $a_X=R_\eta h_X$. Thus, under the isomorphism  def\/ined above the Lie subalgebra $\gc\go\gn(M,\eta)$ of $\gc\go\gn(M,\cald)$ leaving $\eta$ invariant is identif\/ied with the Lie subalgebra $C^\infty(M)^{R_\eta}$ of smooth functions that are invariant under the f\/low of the Reeb vector f\/ield.

Conversely, f\/ixing a contact form $\eta$ the function $h\in C^\infty(M)$ gives a unique {\it Hamiltonian vector field} $X_h\in \gc\go\gn(M,\cald)$ that satisf\/ies $h=\eta(X_h)$. Thus, a {\it Hamiltonian contact structure} is a quadruple $(M,\cald,\eta,h)$. Although any smooth function can be chosen as a Hamiltonian, it is often convenient to choose the function $1=\eta(R_\eta)$ as the Hamiltonian, making the Reeb vector f\/ield $R_\eta$ the Hamiltonian vector f\/ield. We call this a {\it Reeb type Hamiltonian contact structure} and denote it by $(M,\cald,\eta,1)$. It consists only of a contact structure $\cald$ together with a choice of contact form $\eta$ such that $\cald=\ker\eta$ and the Hamiltonian is understood to be the function $1$.

It is sometimes convenient to view $\eta(X)$ in terms of a {\it moment map}. Let $\con(M,\cald)^*$ denote the algebraic dual of $\con(M,\cald)$, and def\/ine the moment map $\Upsilon: \cald^o_+\ra{ 1.6} \con(M,\cald)^*$ by
\begin{gather*}
\langle \Upsilon(p,\eta),X\rangle =\eta(X)(p).
\end{gather*}
So $\Upsilon(p,\eta)\in\con(M,\cald)^*$ is identif\/ied with the linear function ${\rm ev}_p\circ \eta:\con(M,\cald)\ra{1.6} \bbr$ where ${\rm ev}_p$ is the evaluation map at $p$. Fixing the isomorphism $\eta:\con(M,\cald)\ra{1.6} C^\infty(M)$ identif\/ies the image of $\Upsilon$ in $\con(M,\cald)^*$ with the smooth functions $\eta(X)$. We usually consider restricting the moment map to certain f\/inite dimensional Lie subalgebras $\gg$ of $\con(M,\cald)$, and we identify the dual $\gg^*$ with the vector space $\{\eta(X)~|~X\in \gg\}$.

\subsection{First integrals}
We say that a smooth function $f\in C^\infty(M)$ is a {\it first integral} of the contact Hamiltonian structure $(M,\cald,\eta,h)$ if $f$ is constant along the f\/low of the Hamiltonian vector f\/ield, that is if $X_hf=0$. Unlike the symplectic case a contact Hamiltonian is not necessarily a f\/irst integral of its Hamiltonian structure, that is it is not necessarily constant along its own f\/low!

\begin{lemma}\label{1ham}
Let $(M,\cald,\eta,h)$ be a contact Hamiltonian system. Then the following holds:
\begin{gather}\label{1inteqn}
X_hf=(R_\eta h)f+\{h,f\}_\eta.
\end{gather}
In particular, the Hamiltonian function $h$ itself is a first integral if and only if $h\in C^\infty(M)^{R_\eta}$, or equivalently $X_h$ lies in the subalgebra $\gc\go\gn(M,\eta)$.
\end{lemma}

\begin{proof}
We have
\[
X_hf=\pounds_{X_h}\bigl(\eta(X_f)\bigr)=(\pounds_{X_h}\eta)(X_f)+\eta([X_h,X_f])=af+\{h,f\}_\eta
\]
for some $a\in C^\infty(M)$.  But we also have
\[
a=a\eta(R_\eta)=(\pounds_{X_h}\eta)(R_\eta)=(X_h\hook d\eta)(R_\eta)+dh(R_\eta)=R_\eta h
\]
which gives equation (\ref{1inteqn}). For the special case of $f=h$, we have
\begin{gather}\label{Rheqn}
X_hh=(R_\eta h)h.
\end{gather}
So we see that $h$ is constant along its own f\/low if and only if it is constant along the f\/low of the Reeb vector f\/ield. But $R_\eta h=0$ if and only if $h\in C^\infty(M)^{R_\eta}$ which is equivalent to $X_h\in \gc\go\gn(M,\eta)$.
\end{proof}

This begs the question: given a Hamiltonian function $f\in C^\infty(M)$ does there always exist a contact form $\eta$ such that $f$ is constant along the f\/low of the Reeb vector f\/ield of $\eta$, or equivalently given $X\in \con(M,\cald)$ does there always exist an $\eta$ such that $X\in\con(M,\eta)$? The answer is no as now shown.

\begin{proposition}\label{cohprop}
Let $(M,\cald)$ be a co-oriented contact manifold. Then there exist functions $h\in C^\infty(M)$ that are not a first integral of their Hamiltonian system $(M,\cald,\eta,h)$ for any contact form $\eta$.
\end{proposition}

\begin{proof}
At each point of $M$ we know that $R_\eta$ is transversal to $\cald$. So we choose any $h\in C^\infty(M)$ such that at $p\in M$ we have $\ker dh_p=\cald_p$. Then $R_\eta h(p)\neq 0$, and the result follows from Lemma~\ref{1ham}.\footnote{I thank a referee for this very concise proof.}
\end{proof}

In the symplectic case the fact that a Hamiltonian $h$ and a function $f$ commute under Poisson bracket is equivalent to $f$ being a f\/irst integral which is also equivalent to the vector f\/ields $X_h$ and $X_f$ being isotropic with respect to the symplectic form $\gro$. These equivalences no longer hold in the contact case. To have a viable theory in the contact case we should restrict our class of Hamiltonians.

\begin{definition}\label{goodham}
We say that a contact Hamiltonian $h$ is {\it good} if there exists a contact form $\eta$ such that $h$ is constant along the f\/low of the Reeb vector f\/ield $R_\eta$, or equivalently $X_hh=0$. More explicitly, we say that $h$ is a {\it good Hamiltonian with respect to $\eta$}. With this $\eta$ chosen we also say that the contact Hamiltonian system $(M,\cald,\eta,h)$ is {\it good}.
\end{definition}

We now give some straightforward results for good Hamiltonian systems.

\begin{lemma}\label{goodlem}
Let $(M,\cald,\eta,h)$ be a good contact Hamiltonian system. Then
\begin{enumerate}\itemsep=0pt
\item[$1.$] $f\in C^\infty(M)$ is a first integral of $(M,\cald,\eta,h)$ if and only if $\{h,f\}_\eta=0$.
\item[$2.$] If $f$ is a first integral of $h$, then $h$ is a first integral of $f$ if and only if $f$ is constant along the Reeb flow of $R_\eta$, i.e.\ $f$ is a good Hamiltonian with respect to $\eta$.
\item[$3.$] If $f$ is good with respect to $\eta$ then $f$ is a first integral of $h$ if and only if $h$ is a first integral of $f$.
\item[$4.$] Suppose that $h$ and $f$ are mutual first integrals of each other, then $f$ is good with respect to $\eta$.
\item[$5.$] If $f$ is good with respect to $\eta$ then the pointwise linear span of $\{X_h,X_f\}$ is isotropic with respect to $d\eta$ if and only if $\{h,f\}_\eta=0$.
\item[$6.$] If $f$ is a first integral of $(M,\cald,\eta,h)$, then the pointwise linear span of $\{X_h,X_f\}$ is isotropic with respect to $d\eta$ if and only if $X_fh=0$, or equivalently $f$ is good with respect to $\eta$.
\item[$7.$] If $f$ is a first integral of $h$, then $X_h$, $X_f$ span an Abelian subalgebra of $\con(M,\eta)$.
\end{enumerate}
\end{lemma}

\begin{proof}
Items (1)--(4) follow directly from Lemma~\ref{1ham} and Def\/inition~\ref{goodham}. (5) and (6) follow from~(1) and
\begin{gather*}
d\eta(X_h,X_f)=X_hf-X_fh-\{h,f\}_\eta=X_hf-(R_\eta f)h.
\end{gather*}
Finally, Lemma \ref{1ham} implies that $X_h,X_f\in \con(M,\eta)$. Since $f$ is a f\/irst integral of $h$, we have $\eta([X_h,X_f])=\{h,f\}_\eta=0$. But since the only inf\/initesimal contact transformation that is a~section of~$\cald$ is the $0$ vector f\/ield, we have $[X_h,X_f]=0$ which proves~(7).
\end{proof}

\begin{remark}\label{goodistint}
Notice that the Lie algebra $C^\infty(M)^{R_\eta}$ is the set of all good Hamiltonians with respect to $\eta$.
\end{remark}

\begin{remark}\label{abelalg}
More generally the last statement in the proof of Lemma~\ref{goodlem} implies that $\{h,f\}_\eta=0$ if and only if $[X_h,X_f]=0$, and the latter is an Abelian subalgebra of the Lie algebra $\con(M,\cald)$ of inf\/initesimal contact transformations.
\end{remark}

\subsection{Completely integrable contact Hamiltonian systems}

As in symplectic geometry the notion of functions in involution is important; however, in contact geometry it depends on a choice of contact form. For example, given two vector f\/ields $X,Y\in\con(M,\cald)$, a choice of contact form $\eta$ gives a pair of functions $\eta(X),\eta(Y)\in C^\infty(M)$. The vector f\/ields commute if and only if the two functions are in involution. However, choosing a dif\/ferent contact form $\eta'=f\eta$ where $f$ is nowhere vanishing gives two dif\/ferent functions, $f\eta(X)$, $f\eta(Y)$ in involution.

\begin{definition}\label{ininvol} Let $(M,\cald,\eta,h)$ be a contact Hamiltonian system.
A subset $\{h=f_1,f_2,\dots,f_k\}$ of smooth functions is said to be {\it in involution} if $\{f_i,f_j\}_\eta=0$ for all $i,j=1,\dots,k$.
\end{definition}

If $h$ is good with respect to $\eta$ then it follows from (1) of Lemma~\ref{goodlem} that the $f_j$ are all f\/irst integrals of $h$; however, the symmetry of the symplectic case does not hold in general. The function $h$ may not be a f\/irst integral of $f_j$ for $j=2,\dots,k$, since $f_j$ is not necessarily a good Hamiltonian with respect to $\eta$.

\begin{definition}\label{indepfunct}
Let $(M,\cald,\eta,h)$ be a contact Hamiltonian system.
We say that a subset $\{g_1,\dots$, $g_k\}\subset C^\infty(M)$ is {\it independent} if the corresponding set $\{X_{g_1},\dots,X_{g_k}\}$ of Hamiltonian vector f\/ields is pointwise linearly independent on a dense open subset.
\end{definition}

\begin{remark}\label{equivrem}
Unlike the symplectic case this is not equivalent to the condition $dg_1\wedge \cdots\wedge dg_k\neq 0$ on a dense open subset of $M$, since the latter does not hold when one of the Hamiltonian vector f\/ields is $\bbr$-proportional to the Reeb f\/ield whose Hamiltonian is the function $1$. Contact Hamiltonian systems with Hamiltonian equal to $1$ are both interesting and important as we shall see.
\end{remark}

\begin{definition}\label{cidef}
A Hamiltonian contact structure $(M,\cald,\eta,h)$ is said to be {\it completely integrable} if there exists $n+1$ f\/irst integrals, $h,f_1,\dots,f_n$, that are independent and in involution. We denote such a Hamiltonian system by $(M,\cald,\eta,h,\{f_i\}_{i=1}^n)$.
\end{definition}

It follows from equation (\ref{1inteqn}) of Lemma \ref{1ham} that a completely integrable Hamiltonian contact structure $(M,\cald,\eta,h,\{f_i\}_{i=1}^n)$ is automatically good.
However, unlike the symplectic case, $h$~may not be a f\/irst integral of $f_i$, and the subspace spanned by the corresponding vector f\/ields may not be isotropic. From Lemma \ref{goodlem} we have

\begin{proposition}\label{isotropgood}
Suppose that the good Hamiltonian contact structure $(M^{2n+1},\cald,\eta,h)$ has $k+1$ independent first integrals $h=f_0,f_1,\dots,f_k$ with $k\leq n$. Then on a dense open subset the corresponding Hamiltonian vector fields pointwise span a $(k+1)$-dimensional isotropic subspace with respect to $d\eta$ if and only if $f_i$ is good with respect to $\eta$ for all $i=1,\dots,k$.
\end{proposition}

Nevertheless, such contact Hamiltonian structures lift to the usual symplectic Hamiltonian structures on the cone.

\begin{proposition}\label{maxindinv}
Let $(M^{2n+1},\cald,\eta,h)$ be a good contact Hamiltonian system. The maximal number of independent first integrals of the system $(M^{2n+1},\cald,\eta,h)$ that are in involution is $n+1$. In particular, a completely integrable contact Hamiltonian system on $(M,\cald)$ lifts to a~completely integrable Hamiltonian system on $(C(M),\grO)$.
\end{proposition}

\begin{proof}
As discussed at the end of Section~\ref{congroupsect} we can lift the independent Hamiltonian vector f\/ields of the f\/irst integrals to the symplectization $(C(M),\grO)$, and a direct computation~\cite{Boy10a} shows that an Abelian subalgebra of $\con(M,\cald)$ lifts to an Abelian subalgebra in $\gs\gy\gm(C(M),\grO)$. The maximal dimension of such an Abelian subalgebra is $n+1$.
\end{proof}

Generally, the Hamiltonian vector f\/ields may not be complete, so they do not integrate to an element of the contactomorphism group $\Con(M,\cald)$, but only to the corresponding pseudogroup. Even if the manifold is compact so the Abelian Lie algebra of Hamiltonian vector f\/ields integrates to an Abelian group $\cala\subset \Con(M,\cald)$, it may not be a closed Lie subgroup of $\Con(M,\cald)$ as we shall see in Example \ref{cosphere} below.

\begin{definition}\label{goodcomint}
We say that the completely integrable contact Hamiltonian system $(M,\cald,\eta,h$, $\{f_i\}_{i=1}^n)$ is {\it completely good} if the f\/irst integral $f_i$ is good with respect to $\eta$ for all $i=1,\dots,n$.
\end{definition}

From the def\/initions it follows that
\begin{proposition}\label{iffcompgood}
A completely integrable contact Hamiltonian system $(M,\cald,\eta,h=f_0,f_i,\dots, f_n)$ is completely good if and only if the corresponding commuting Hamiltonian vector fields $X_{f_0},\dots$, $X_{f_n}$ lie in the subalgebra $\con(M,\eta)$.
\end{proposition}

Most of the known completely integrable contact Hamiltonian systems are completely good, but here is a simple example which is not completely good.

\begin{example}\label{notgood}
Take $M=\bbr^3$ with standard coordinates $(x,y,z)$ and standard contact form $\eta=dz-ydx$. The Reeb vector f\/ield is $R_\eta=\d_z$. Take $h=-y$ as the Hamiltonian which is good with respect to $\eta$ and take $f=z$ as a f\/irst integral. The f\/irst integral $f$ is not good with respect to $\eta$ since $R_\eta f=\d_zf=1$. The functions $h$ and $f$ are in involution, since the corresponding vector f\/ields, which are
\[
X_h=\d_x, \qquad X_f=z\d_z+y\d_y,
\]
commute. However, they are not isotropic with respect to $d\eta$ since $d\eta(X_h,X_f)=y=-h\neq 0$.

A question that arises is whether there exists a dif\/ferent contact form $\eta'$ in $\cald$ such that the system $(\bbr^3,\cald,\eta',h,f)$ is completely good. So we look for a smooth positive function $g$ on $\bbr^3$ such that $\eta'=g\eta$ that satisf\/ies both $\pounds_{X_f}(g\eta)=0$ and $\pounds_{X_h}(g\eta)=0$. This implies $g_x=0$ and
\[
0=\pounds_{X_f}(g\eta)=(X_fg)\eta+g\pounds_{X_f}\eta=\bigl(X_fg+g\bigr)\eta.
\]
Thus, $g$ must be independent of $x$ and homogeneous of degree $-1$ in $y$ and $z$. But there is no such positive smooth function on $\bbr^3$. So this completely integrable contact Hamiltonian system is not completely good with respect to any contact form representing $\cald=\ker \eta$. Notice also that the same argument shows that the Hamiltonian $z$ is not a good Hamiltonian with respect to any $\eta$ representing the standard contact structure $(M,\cald=\ker \eta)$.

Nevertheless this contact Hamiltonian system lifts to a completely integrable Hamiltonian system on the symplectic cone $C(M)\approx \bbr^3\times \bbr^+$. In coordinates $(x,y,z,r)$ the symplectic form is
\[
\grO=r^2dx\wedge dy +2rdr\wedge (dz-ydx).
\]
The lifted vector f\/ields are $\hX_h=\d_x$ and $\hX_f=z\d_z+y\d_y-\frac{1}{2}r\d_r$ with Hamiltonians~$-r^2y$ and~$r^2z$, respectively.
\end{example}

\begin{definition}\label{comintReeb}
A completely integrable contact Hamiltonian system $(M^{2n+1},\cald,\eta,f_0=h,f_1,\dots$, $f_n)$ is of {\it Reeb type} if $f_i=1$ for some $i=0,\dots,n$.
\end{definition}

Of course, this is equivalent to the condition that the Reeb vector f\/ield $R_\eta$ lies in the Abelian Lie algebra spanned by the Hamiltonian vector f\/ields $X_{f_0},\dots,X_{f_n}$.
We have
\begin{theorem}\label{reebcomgood}
A completely integrable contact Hamiltonian system $(M^{2n+1},\cald,\eta,f_0=h,f_1,\dots,$ $f_n)$ of Reeb type is completely good.
\end{theorem}

\begin{proof}
It is well known that the centralizer of the Reeb vector f\/ield $R_\eta$ of a contact form $\eta$ in $\con(M,\cald)$ is the subalgebra $\con(M,\eta)$. So the condition $\{1,f_i\}_\eta=0$ implies that $X_i\in\con(M,\eta)$ for all $i$ and this implies $R_\eta f_i=0$ for all $i$.
\end{proof}

The converse does not hold. Example \ref{cosphere} below is a completely good contact Hamiltonian system on a compact manifold that is not of Reeb type. The contact Hamiltonian systems that we are mainly concerned with in this paper are completely good.

\begin{definition}\label{torictype}
A completely integrable contact Hamiltonian system $(M^{2n+1},\cald,\eta,f_0,f_1,\dots,f_n)$ is said to be of {\it toric type} if the corresponding Hamiltonian vector f\/ields $X_{f_0},\dots,X_{f_n}$ form the Lie algebra of a torus $T^{n+1}\subset \Con(M,\cald)$. In this case we also call $(M^{2n+1},\cald,\eta)$ a {\it toric contact manifold}.
\end{definition}

\begin{example}\label{cosphere}
Consider the unit sphere bundle $S(T^*T^{n+1})$ of the cotangent bundle of an $(n+1)$-torus $T^{n+1}$. In the canonical coordinates $(x^0,\dots,x^n;p_0,\dots,p_n)$ on the cotangent bundle, $S(T^*T^{n+1})$ is represented by $\sum_{i=0}^np_i^2=1$, with $(x^0,\dots,x^n)$ being the coordinates on the torus $T^{n+1}$. It is easy to see that the restriction of the canonical 1-form $\theta=\sum_ip_idx^i$ on $T^*T^{n+1}$ to $S(T^*T^{n+1})$ is a contact form $\eta=\theta|_{S(T^*T^{n+1})}$ on $S(T^*T^{n+1})$.
Moreover, this contact structure is toric since $T^{n+1}$ acts freely on $S(T^*T^{n+1})$ and leaves $\eta$ invariant. The Reeb vector f\/ield~$R_\eta$ of $\eta$ is the restriction of $\sum_{i=0}^np_i\d_{x^i}$ to $S(T^*T^{n+1})$, and this does not lie in the Lie algebra~$\gt_{n+1}$ of the torus which is spanned by $\d_{x^i}$. So this toric contact structure is not of Reeb type. Note that the vector f\/ields~$\{R_\eta,\d_{x^0},\dots,\d_{x^n}\}$ form an $(n+2)$-dimensional Abelian Lie algebra, but the independence condition of Def\/inition~\ref{indepfunct} fails. Note, however, that the vector f\/ields $\{R_\eta,\d_{x^1},\dots,\d_{x^n}\}$ do form an $(n+1)$-dimensional Abelian Lie algebra and they are independent on the dense open subset $(p_0\neq 0)$ in $S(T^*T^{n+1})$. This gives a completely good integrable system of Reeb type on~$S(T^*T^{n+1})$ which, however, is not of toric type, since the vector f\/ield~$R_\eta$ generates an~$\bbr$ action, and there are orbits of this $\bbr$ action whose closure is $T^{n+1}$. So the subgroup $T^n\times \bbr$ generated by this completely integrable system is not a closed Lie subgroup of $\Con(S(T^*T^{n+1}),\eta)$.
\end{example}

\begin{example}\label{Heis}
An example of a completely good integrable contact Hamiltonian system of Reeb type that is not toric is given by the standard Sasakian contact structure on the Heisenberg group $\gH^{2n+1}$ \cite{Boy09}. As a contact manifold $\gH^{2n+1}$ is just $\bbr^{2n+1}$ with contact form
\[
\eta=dz-\sum_{i=1}^ny_idx^i
\]
given in global coordinates $(x^1,\dots,x^n,y_1,\dots, y_n,z)$. The connected component of the Sasakian automorphism group is the semi-direct product $(U(n)\times \bbr^+)\ltimes \gH^{2n+1}.$ Let $X_i$ denote the vector f\/ields corresponding to the diagonal elements of $U(n)$. Then the functions $(1,\eta(X_1),\dots,\eta(X_n))$ make $(\gH^{2n+1},\cald,\eta,1)$ into a completely integrable contact Hamiltonian system of Reeb type which is not toric since the corresponding Abelian group is $T^n\times \bbr$ where the Reeb vector f\/ield $R_\eta=\d_z$ generates a real line $\bbr$.

Another completely integrable contact Hamiltonian system on $\gH^{2n+1}$ is obtained by taking the Hamiltonian as a linear combination of the functions $(\eta(X_1),\dots,\eta(X_n))$, say the sum $h=\sum_i\eta(X_i)$, and adding the function $\eta(D)$ where $D$ is the generator of the $\bbr^+$ gives a completely integrable contact Hamiltonian system on $\gH^{2n+1}$ which by Proposition~\ref{iffcompgood} is not completely good, since $D\not\in \con(M,\eta)$.
\end{example}

Completely integrable contact structures of toric type were f\/irst studied in~\cite{BM93}, but their def\/inition was local in nature, hence more general than ours, involving the condition of transverse ellipticity. The passage from local to global involves the vanishing of the monodromy of the so-called Legendre lattice. On compact manifolds global toric actions were classif\/ied by Lerman in~\cite{Ler02a} where it is shown that the toric contact structures of Reeb type are precisely those having a description in terms of polyhedra as in  the symplectic case.  By averaging over the torus one constructs a contact form in the contact structure $\cald$ that is invariant under the action of the torus. Thus,

\begin{proposition}\label{toricgood}
For a completely integrable contact Hamiltonian system of toric type there is a~contact form $\eta$ such that the completely integrable contact Hamiltonian system $(M^{2n+1},\cald,\eta,f_0,$ $f_1,\dots,f_n)$ is completely good.
\end{proposition}

More generally one can use a slice theorem \cite{Ler02a,BG05} to prove
\begin{proposition}\label{propergood}
Let $(M^{2n+1},\cald,\eta,f_0,f_1,\dots,f_n)$ be a completely integrable contact Hamiltonian system, and suppose that the corresponding Hamiltonian vector fields $\{X_{f_i}\}_{i=0}^n$ form the Lie algebra of an $(n+1)$-dimensional Abelian subgroup $\cala$ of $\Con(M,\cald)$ whose action on $M$ is proper. Then there exists a contact form $\eta_0$ representing $\cald$ such that $\cala\subset \Con(M,\eta_0)$ and the corresponding completely integrable contact Hamiltonian system  $(M^{2n+1},\cald,\eta_0,f_0,f_1,\dots,f_n)$ is completely good.
\end{proposition}

The Lie algebra $\ga$ of $\cala$ provides $M$ with a singular foliation $\calf_\ga$ whose generic leaves are submanifolds of maximal dimension $n+1$. Note that $\calf_\ga$ is a true $(n+1)$-dimensional foliation on a dense open subset $W\subset M$. Let $\pi_\cald:TM\ra{1.7} \cald$ denote the natural projection. Then the image $\pi_\cald(\ga)$ def\/ines another singular foliation on $M$ called the {\it Legendre foliation} whose generic leaves are $n$-dimensional, called {\it Legendrian submanifolds}. The leaves of both of these foliations are endowed with a canonical af\/f\/ine structure \cite{Lib91,Pan90}. In particular, compact leaves are tori. This is part of the contact analogue of the well known Arnold--Liouville theorem in symplectic geometry.

\subsection{Conjugacy classes of contact Hamiltonian systems}\label{conjclass}
Let $(M^{2n+1},\cald,\eta,f_0,f_1,\dots,f_n)$ be a completely integrable contact Hamiltonian system. Then the vector f\/ields $\{X_{f_i}\}_{i=0}^n$ span an $(n+1)$-dimensional Abelian Lie subalgebra $\ga_{n+1}$ of $\con(M,\cald)$. Such Lie subalgebras are maximal in the sense that they have the maximal possible dimension for Abelian subalgebras of $\con(M,\cald)$. We are interested in the conjugacy classes of completely integrable Hamiltonian systems. These should arise from conjugacy classes of maximal Abelian subalgebras of $\con(M,\cald)$. First notice that by indentifying $\con(M,\cald)$ as the Lie algebra of left invariant (or right invariant) vector f\/ields on $M$ leaving $\cald$ invariant, we can make the identif\/ication ${\rm Ad}_\phi X=\phi_*X$ for $X\in \con(M,\cald)$ and $\phi\in\Con(M,\cald)$. Thus, by duality we have ${\rm Ad}^*_\phi\eta=(\phi^{-1})^*\eta$. So under conjugation if $h=\eta(X)$ and $h'=((\phi^{-1})^*\eta)(\phi_*X)$ we have
\begin{gather*}
h'(p) = ((\phi^{-1})^*\eta)(\phi_*X)(p)= \eta(\phi^{-1}_*\phi_*X)(\phi^{-1}(p))= \eta(X)(\phi^{-1}(p))
  = (h\circ \phi^{-1})(p).
\end{gather*}
This leads to

\begin{definition}\label{equivcompint}
We say the contact Hamiltonian systems $(M,\cald,\eta,h)$ and $(M,\cald,\eta',h')$ are {\it conjugate} if there  exists $\phi\in\Con(M,\cald)$ such that $\eta'=(\phi^{-1})^*\eta$ and $h'=h\circ \phi^{-1}$.
\end{definition}

Now goodness is preserved under conjugation. More explicitly,

\begin{lemma}\label{goodconj}
If $h\in C^\infty(M)$ is good with respect to the contact form $\eta$, then $h'=h\circ \phi^{-1}$ is good with respect to $\eta'=(\phi^{-1})^*\eta$ for any $\phi\in\Con(M,\cald)$.
\end{lemma}

\begin{proof}
From equation (\ref{Rheqn}) any smooth function $h$ is good with respect to $\eta$ if and only if $X_hh=0$. We have
\begin{gather*}
X_{h'}h' = ({\rm Ad}_\phi X_h)(h\circ \phi^{-1})=(\phi_*X_h)(h\circ \phi^{-1})= d((\phi^{-1})^*h) (\phi_*X_h) \\
\hphantom{X_{h'}h'}{}  = ((\phi^{-1})^*dh)(\phi_*X_h)=dh(\phi^{-1}_*\phi_*X_h)=dh(X_h)=X_hh=0.\tag*{\qed}
\end{gather*}\renewcommand{\qed}{}
\end{proof}

Similarly, two completely integrable contact Hamiltonian systems $(M,\cald,\eta,h,\{f_i\}_{i=1}^n)$ and $(M,\cald,\eta',h',\{f'_i\}_{i=1}^n)$ are {\it conjugate} if there exists $\phi\in\Con(M,\cald)$ such that $\eta'=(\phi^{-1})^*\eta$, $h'=h\circ \phi^{-1}$, and $f'_i=f_i\circ \phi^{-1}$ for all $i=1,\dots,n$. It is easy to see that

\begin{lemma}\label{goodconj2}
If a completely integrable contact Hamiltonian system $(M^{2n+1},\cald,\eta,f_0,f_1,\dots,f_n)$ is either completely good, of Reeb type, or toric type, so is any conjugate of $(M^{2n+1},\cald,\eta,f_0,f_1,\dots,$ $f_n)$ by any element of $\Con(M,\cald)$.
\end{lemma}

In the case of a completely good integrable contact Hamiltonian system we can consider the conjugacy of an entire Abelian algebra instead of the functions individually. In this case the conjugate system $(M^{2n+1},\cald,\eta',f'_0,f'_1,\dots,f'_n)$ satisf\/ies{\samepage
\begin{gather*}
\eta'=\big(\phi^{-1}\big)^*\eta, \qquad f'_i=\sum_{j=0}^nA_{ij}f_j\circ\phi^{-1} ~\text{for all $i=0,\dots,n$,}
\end{gather*}
where $A_{ij}$ are the components of a matrix $A\in {\rm GL}(n+1,\bbr)$.}

Conjugacy classes of completely good integrable contact Hamiltonian systems correspond to conjugacy classes of maximal Abelian subalgebras of dimension $n+1$ in $\con(M,\cald)$ whose projections onto $\cald$ are isotropic with respect to $d\eta$ for any (hence, all) contact form $\eta$ representing~$\cald$. The completely integrable contact Hamiltonian systems considered in Section~\ref{torics2s3} are toric, hence completely good. In particular, we are interested in the inequivalent ways that a given contact Hamiltonian system $(M^{2n+1},\cald,\eta,h)$ may be completely integrable. This corresponds to maximal Abelian subalgebras of $\con(M,\cald)$ containing a f\/ixed one-dimensional subalgebra generated by~$X_h$. We give an important example of this in Section~\ref{torics2s3} where the distinct conjugacy classes are related to inequivalent Sasaki--Einstein metrics on $S^2\times S^3$. In this case the Hamiltonian vector f\/ield is the Reeb vector f\/ield of a preferred contact form, and the contact Hamiltonian systems are actually $T^2$-equivariantly equivalent, but not $T^3$-equivariantly equivalent.

The general situation we are interested in is as follows: consider a contact structure $\cald$ on a~manifold $M$ of dimension $2n+1$, and f\/ix a contact form $\eta$. Suppose that the group $\Con(M,\eta)$ has exactly $\gn(\cald,n+1)$ conjugacy classes of maximal tori of dimension $n+1$. Thus, the Hamiltonian system $(M,\cald,\eta,1)$ has at least $\gn(\cald,n+1)$ inequivalent ways that make it completely integrable. Moreover, it is easy to construct examples \cite{Boy10b,Boy10a} where $\gn(\cald,n+1)$ is quite large. In fact, for the toric contact structures $Y^{p,q}$ on $S^2\times S^3$ discussed in Section~\ref{torics2s3}, there are at least $p-1$ conjugacy classes of maximal tori of dimension $3$ when $p$ is prime. So the corresponding f\/irst integrals cannot all be independent. The same is true on the level of the symplectic base space discussed at the end of Section~\ref{section1.1}. A question that arose during the conference is whether this phenomenon is at all related to that of superintegrability. However, it appears that this is not the case. In symplectic geometry superintegrable Hamiltonian systems arose from the noncommutative or generalized Arnold--Liouville theorem of Nehorosev~\cite{Neh72}, and Mishchenko and Fomenko~\cite{MiFo78}, and recently this notion has come to the forefront of mathematical physics (see for example~\cite{TWHMPR04,Fas05} and references therein). Nevertheless, it would be of interest to develop the non-commutative theory in the contact setting.

\section[Toric contact structures on $S^2\times S^3$]{Toric contact structures on $\boldsymbol{S^2\times S^3}$}\label{torics2s3}

In this section we consider the special case of toric contact structures, denoted by $Y^{p,q}$ where $p$, $q$ are relatively prime integers satisfying $0<q<p$, on $S^2\times S^3$ that were f\/irst constructed in \cite{GMSW04a}. The contact bundle of these structures has vanishing f\/irst Chern class; however, they are not the most general toric contact structures with vanishing f\/irst Chern class. The latter were studied in \cite{CLPP05,MaSp05b}. Here for reasons of simplicity we consider only the case of the $Y^{p,q}$. More general toric contact structures, including ones with nonvanishing f\/irst Chern class, on $S^3$ bundles over $S^2$ are considered in greater detail in~\cite{BoPa10}.

\subsection[Circle reduction of $S^7$]{Circle reduction of $\boldsymbol{S^7}$}\label{circred}

The structures $Y^{p,q}$ on $S^2\times S^3$ can be obtained by the method of symmetry reduction from the standard contact structure on~$S^7$ by a certain circle action. For a review of contact reduction we refer to Chapter~8 of~\cite{BG05}, while complete details of this case can be found in~\cite{BoPa10}. Consider the standard~$T^4$ action on~$\bbc^4$ given by~$z_j\mapsto e^{i\theta_j}z_j$. Its moment map $\Upsilon_4:\bbc^4\setminus \{0\}\ra{1.7} \gt^*_4=\bbr^4$
is given by
\begin{gather*}
\Upsilon_4(z)=\big(|z_1|^2,|z_2|^2,|z_3|^2,|z_4|^2\big).
\end{gather*}
Now we consider the circle group $T(p,q)$ acting on $\bbc^4\setminus \{0\}$ by
\begin{gather}\label{s1action}
(z_1,z_2,z_3,z_4)\mapsto \big(e^{i(p-q)\theta}z_1,e^{i(p+q)\theta}z_2,e^{-ip\theta}z_3, e^{-ip\theta}z_4\big),
\end{gather}
where $p$ and $q$ are positive integers satisfying $1\leq q<p$. The moment map for this circle action is given by
\begin{gather*}
\Upsilon_1(z)=(p-q)|z_1|^2+(p+q)|z_2|^2 -p\big(|z_3|^2+|z_4|^2\big).
\end{gather*}
Representing $S^7$ by $(p-q)|z_1|^2+(p+q)|z_2|^2+p(|z_3|^2+|z_4|^2)=1$ shows that the zero set of $\Upsilon_1$ restricted to $S^7$ is $S^3\times S^3$ represented by
\begin{gather*}
(p-q)|z_1|^2+(p+q)|z_2|^2=\frac{1}{2}, \qquad |z_3|^2+|z_4|^2=\frac{1}{2p}.
\end{gather*}
The action of the circle $T(p,q)$ on this zero set is free if and only if $\gcd(q,p)=1$. Assuming  this the procedure of contact reduction gives the quotient manifold $Y^{p,q}=S^7/T(p,q)$ with its induced contact structure $\cald_{p,q}$.

\begin{proposition}\label{ypq}
The quotient contact manifold $(Y^{p,q},\cald_{p,q})$ is diffeomorphic to $S^2\times S^3$.
\end{proposition}

\begin{proof}({\it Outline.})
By a result of Lerman \cite{Ler04} the condition $\gcd(q,p)=1$ implies that the manifold $Y^{p,q}$ is simply connected and that $H_2(Y^{p,q},\bbz)  =\bbz$. Furthermore, it is not dif\/f\/icult to see \cite{BoPa10} that the f\/irst Chern class $c_1(\cald_{p,q})$ vanishes, so the manifold is spin. The result then follows by the Smale--Barden classif\/ication of simply connected 5-manifolds.
\end{proof}

It follows from \cite{MaSp06} that the contact structures $Y^{p,q}$ are precisely the ones discovered in~\cite{GMSW04a} that admit Sasaki--Einstein metrics. We should also mention that the case $(p,q)=(1,0)$ admits the well-known homogeneous Sasaki--Einstein metric found over 30 years ago by Tanno~\cite{Tan79}. Our proof expresses the $Y^{p,q}$ as $S^1$-orbibundles over certain orbifold Hirzebruch surfaces. A~connection between the $Y^{p,q}$ and Hirzebruch surfaces was anticipated by Abreu \cite{Abr09}.

\subsection{Outline of the proof of Theorem~\ref{SEequiv}}\label{proofthm}

The proof that the contact structures are inequivalent if $p'\neq p$ uses contact homology as very brief\/ly described in Section~\ref{conhom}. For full details of a more general result we refer to~\cite{BoPa10}; however, the case for the $Y^{p,q}$ was worked out by Abreu and Macarini~\cite{AbMa10}, so we shall just refer to their paper for the proof.

To prove that the contact structures $Y^{p,q}$ and $Y^{p,q'}$ are contactomorphic requires a judicious choice of Reeb vector f\/ield in the Sasaki cone, or equivalently, a judicious choice of contact form. First the inf\/initesimal generator of the circle action given by equation~(\ref{s1action}) is $L_{p,q}=(p-q)H_1+(p+q)H_2-p(H_3+H_4)$ where $H_j$ is the inf\/initesimal generator of the action $z_j\mapsto e^{i\theta_j}z_j$ on $\bbc^4$ restricted to $S^3\times S^3$. Choosing the vector f\/ield $R_{p,q}=(p+q)H_1+(p-q)H_2+p(H_3+H_4)$ gives the $T^2$-action generated by $L_{p,q}$ and $R_{p,q}$ as
\begin{gather*}\label{t2act}
\bfz\mapsto \big(e^{i((p+q)\phi+(p-q)\theta)}z_1,e^{i((p-q)\phi+(p+q)\theta}z_2,e^{ip(\phi-\theta)}z_3, e^{ip(\phi-\theta)}z_4\big).
\end{gather*}
Upon making the substitutions $\psi=\phi-\theta$ and $\chi=(p-q)\psi+2p\theta$ we obtain the action
\begin{gather}\label{t2genact1}
\bfz\mapsto \big(e^{i(2q\psi+\chi)}z_1,e^{i\chi}z_2,e^{ip\psi}z_3, e^{ip\psi}z_4\big).
\end{gather}
It is easy to check that $R_{p,q}$ is a Reeb vector f\/ield, that is $\eta(R_{p,q})>0$ everywhere, where $\eta$ is the contact form on~$Y^{p,q}$ induced by the reduction procedure. Furthermore, the circle action on~$Y^{p,q}$ generated by~$R_{p,q}$ is quasiregular and gives an orbifold Boothby--Wang quotient space $\calz_{p,q}$. We now identify this quotient as an orbifold Hirzebruch surface, that is a Hirzebruch surface with a nontrivial orbifold structure. To do so we equate the $T^2$ reduction of $S^3\times S^3$ by the action (\ref{t2genact1}) with the $\bbc^*\times\bbc^*$ quotient of $\bbc^2\setminus \{0\}\times \bbc^2\setminus \{0\}$ by the complexif\/ied $T_\bbc^2$-action with parameters $\grt$, $\grz$:
\begin{gather}\label{C*act}
\bfz\mapsto \big(\grt^{2q}\grz z_1,\zeta z_2,\grt^{p}z_3, \grt^{p}z_4\big).
\end{gather}
Since $p$ and $q$ are relatively prime there are two cases to consider. If $p$ is odd then equation (\ref{C*act}) stays as is; however, if $p$ is even we can redef\/ine $\grt$ to give
\begin{gather*}
\bfz\mapsto \big(\grt^{q}\grz z_1,\zeta z_2,\grt^{\tp}z_3, \grt^{\tp}z_4\big),
\end{gather*}
where $p=2\tp$.

When $p$ is odd we def\/ine homogeneous coordinates in $\bbc\bbp^1\times \bbc\bbp^2$ by setting $(w_1,w_2)=(z_3,z_4)$ and $(y_1,y_2,y_3)=(z_2^{p}z_3^{2q},z_2^{p}z_4^{2q},z_1^{p})$, we see that the quotient space $\calz_{p,q}$ is represented by the equation
\begin{gather*}
w_1^{2q}y_2=w_2^{2q}y_1.
\end{gather*}
This is well known to be the even Hirzebruch surface $S_{2q}$ which is dif\/feomorphic to $S^2\times S^2$. However, the divisors $E=(z_1=0)$ and $F=(z_2=0)$ are branch divisors with ramif\/ication index $p$. So we have a non-trivial orbifold structure $(S_{2q},\grD_p)$ with
\begin{gather*}
\grD_p=\left(1-\frac{1}{p}\right)(E+F).
\end{gather*}
When $p$ is even we get an orbifold Hirzebruch surface $(S_q,\grD_{\tp})$ with
\begin{gather*}
\grD_{\tp}=\left(1-\frac{1}{\tp}\right)(E+F).
\end{gather*}
In this case $q$ must be odd so we have an odd Hirzebruch surface $S_q$ which is dif\/feomorphic to the blow-up of $\bbc\bbp^2$ at a point.

\begin{remark}\label{logdelPez}
The orbifold structure plays a crucial role in the existence of Sasaki--Einstein metrics in the Sasaki cone of $Y^{p,q}$. In both cases the surfaces $(S_{2q},\grD_{p})$ and $(S_q,\grD_{\tp})$ are log del Pezzo, that is, the orbifold anticanonical divisor is ample; whereas, it is well known that the ordinary anticanonical divisor for a Hirzebruch surface $S_k$ is ample only for $k=0,1$. So the orbifold structure allows for positive Ricci curvature metrics on $Y^{p,q}$ realized as the total space of the $S^1$-orbibundle or alternatively Seifert f\/ibration over the log del Pezzo surfaces $(S_k,\grD_p)$.
\end{remark}

Returning to the proof of Theorem~\ref{SEequiv}, we generalize a result of Karshon \cite{Kar99,Kar03} to show that $Y^{p,q}$ and $Y^{p,q'}$ are contactormophic. Let us assume that $p$ is odd as the proof in the $p$ even case is similar. Karshon \cite{Kar03} shows that the smooth manifold $S^2\times S^2$ with the symplectic form $\gro=\gro_1+k\gro_2$ where $\gro_i$ is the standard K\"ahler form on the ith copy of $S^2$ has exactly $k$ conjugacy classes of $2$-tori in its symplectomorphism group, and that these correspond to dif\/ferent choices of compatible complex structures. In this way each even Hirzebruch surface $S_{2q}$ occurs with $0\leq q< k$, and using her previous results \cite{Kar99} Karshon shows that they are all $S^1$ equivariantly symplectomorphic.

In our case we have the orbifold Hirzebruch surface $(S_{2q},\grD_p)$ with the same symplectic form $\gro=\gro_1+p\gro_2$ as in the smooth case. Now let us write $(S_{2q},\varnothing)$ to denote the even Hirzebruch surface with the trivial orbifold structure, that is the charts are just the standard manifold charts. In this situation as emphasized in \cite{GhKo05} we consider the map $\BOne:(S_{2q},\grD)\ra{1.6} (S_{2q},\varnothing)$ which is the identity as a set map, and a Galois cover with trivial Galois group. Now let $1\leq q,q'<p$ and $K:S_{2q}\ra{1.6} S_{2q'}$ be the Karshon symplectomorphism. Then we have a commutative diagram
\begin{gather*}
\begin{matrix}(S_{2q},\grD_{p}) &\fract{K_p}{\ra{2.8}} &(S_{2q'},\grD_{p}) \\
                         \decdnar{\BOne} &&           \decupar{\BOne^{-1}} \\
                         (S_{2q},\varnothing)&\fract{K}{\ra{2.8}} &(S_{2q'},\varnothing)
\end{matrix}
\end{gather*}
which def\/ines the upper horizontal arrow $K_p$ and shows that it too is an $S^1$-equivariant symplectomorphism. But it is easy to see that $K$ leaves both divisors $E$ and $F$ invariant separately, and this shows that $K_p$ is an orbifold dif\/feomorphism, and hence, an orbifold symplectomorphism. But then as shown in \cite{Ler03b,Boy10a} this symplectomorphism lifts to a $T^2$-equivariant contactomorphism. Hence, the contact structures $Y^{p,q}$ and $Y^{p,q'}$ are contactomorphic. This completes the proof of Theorem~\ref{SEequiv}.

\subsection[Extremal Sasakian metrics on $S^2\times S^3$]{Extremal Sasakian metrics on $\boldsymbol{S^2\times S^3}$}
Extremal Sasakian structures are closely related to the better known extremal K\"ahlerian structures \cite{Cal82}, and were f\/irst def\/ined in \cite{BGS06}. They were studied further in \cite{BGS07b,Boy10b,Leg10}.

Consider the following deformation of a Sasakian structure $\cals=(\xi,\eta,\Phi,g)$ by deforming the contact form by sending $\eta\mapsto \eta(t)=\eta+t\grz$ where $\grz$ is a basic 1-form with respect to the characteristic foliation $\calf_\xi$ def\/ined by the Reeb vector f\/ield $R_\eta=\xi$, that is, $\grz$ depends only on the variables transverse to the Reeb f\/ield $\xi$. Note that since the $\grz$ is basic $\eta(t)$ has the same Reeb vector f\/ield $\xi$ as $\eta$ for all $t$ in a suitable interval containing $0$ and such that $\eta(t)\wedge d\eta(t)\neq 0$. Let $L_\xi$ denote the 1-dimensional subbundle of the tangent bundle generated by $\xi$, and set $\nu(\calf_\xi)=TM/L_\xi$. Let $\bJ$ denote the complex structure on $\nu(\calf_\xi)$ induced by $\Phi$, and let ${\mathcal S}(\xi, \bar{J})$ denote the deformation space of such Sasakian structures. Following \cite{BGS06} we let $s_g$ denote the scalar curvature of $g$ and def\/ine
the ``energy functional'' $E:{\mathcal S}(\xi,\bar{J})\ra{1.4} \bbr$ by
\begin{gather}\label{var}
E(\cals) =  \int _M s_g ^2 d{\mu}_g ,
\end{gather}
i.e. the $L^2$-norm of the scalar curvature. Critical points $\cals$ of this functional are called {\it extremal Sasakian structures}. Similar to the K\"ahlerian case, the Euler-Lagrange equations for this functional give \cite{BGS06} the following statement. A Sasakian structure $\cals\in {\mathcal S}(\xi,\bar{J})$ is a critical point for the energy functional (\ref{var}) if and only if the gradient vector f\/ield $\partial^\#_gs_g$ is transversely holomorphic. In particular, Sasakian metrics with constant scalar curvature are extremal. Here $\partial^\#_g$ is the $(1,0)$-gradient vector f\/ield def\/ined by $g(\partial^\#_g\varphi,\cdot)= \bar{\partial}\varphi$.
It is important to note that a Sasakian metric $g$ is extremal if and only if the `transverse metric' $g_\cald=d\eta\circ (\Phi\otimes \BOne)$ is extremal in the K\"ahlerian sense which follows from the well known relation (cf. \cite{BG05}) between scalar curvatures $s_g=s_g^T-2n$ where $s_g^T$ is the scalar curvature of the transverse metric $g_\cald$.

We now consider the toric contact structures $Y^{p,q}$ on $S^2\times S^3$. As mentioned in the introduction it is known \cite{GMSW04a} that for every relatively prime pair $(p,q)$ with $1\leq q<p$ there is a unique Sasaki--Einstein metric in the Sasaki cone of $Y^{p,q}$. Sasaki--Einstein metrics have constant scalar curvature, and thus, def\/ine extremal Sasakian structures. In fact when $c_1(\cald)=0$ any Sasakian metric of constant scalar curvature is related by a transverse homothety to a Sasaki--Einstein metric. So there is a unique ray of Sasakian metrics with constant scalar curvature. Now f\/ixing the contact structure $\cald_p$ up to isomorphism there are $\phi(p)$ inequivalent toric contact structures associated with $\cald_p$, hence, $\phi(p)$ Sasaki cones $\gt^+_3(q)$, and each Sasaki cone has a unique ray of constant scalar curvature metrics. Moreover, by applying the Openness theorem of \cite{BGS06} there is an open cone $\ge_3(q)\subset \gt^+_3(q)$ of extremal Sasaki metrics about each of these rays.

Other extremal Sasakian metrics in the toric contact structures $Y^{p,q}$ can be obtained by a~dif\/ferent choice of Reeb vector f\/ield.

\begin{proposition}\label{otherext}
Consider the toric contact structure $(Y^{p,q},\cald_p)$ on $S^2\times S^3$ obtained from the circle reduction of Section~{\rm \ref{circred}}. Let $\eta$ be the contact form whose Reeb vector field is
\[
R_\eta =(p-q)H_1+(p+q)H_2+p(H_3+H_4).
\]
Then the induced Sasakian structure $\cals=(R_\eta,\eta,\Phi,g)$ is extremal with non-constant scalar curvature.
\end{proposition}

\begin{proof}
The quotient space of $Y^{p,q}$ by the two torus $T^2$ generated by the vector f\/ields $L_{p,q}=(p-q)H_1+(p+q)H_2-p(H_3+H_4)$ and $R_\eta =(p-q)H_1+(p+q)H_2+p(H_3+H_4)$ is
\[
\bbc\bbp(\bp_-,\bp_+)\times \bbc\bbp^1
\]
with the product complex structure where $\bbc\bbp(\bp_-,\bp_+)$ is a weighted projective space, and $(\bp_-,\bp_+)=(p-q,p+q)$ if $p-q$ is odd, and $(\bp_-,\bp_+)=(\frac{p-q}{2},\frac{p+q}{2})$ if $p-q$ is even. The orbifold K\"ahler form on $\bbc\bbp(\bp_-,\bp_+)\times \bbc\bbp^1$ is
\begin{gather}\label{prodKah}
\gro=p\gro_1+\gcd(p-q,p+q)\gro_2,
\end{gather}
where $\gro_i$ is the standard Fubini--Study form on the $i$th factor. Note that because of the conditions on $(p,q)$ the $\gcd(p-q,p+q)$ is either $1$ or $2$, and it is $2$ only when $p$ and $q$ are both odd.

Now according to Bryant \cite{Bry01} weighted projective spaces have Bochner-f\/lat K\"ahler orbifold metrics which by \cite{DaGa06} are extremal. Moreover, when the weights are not all one the extremal metrics have non-constant scalar curvature. The standard K\"ahler structure $\gro_2$ on $\bbc\bbp^1$ has constant scalar curvature, so the weighted product structure $\gro$ of equation (\ref{prodKah}) has the form $s_\gro=p^{-1}s_1+c$ where $s_1$ is the scalar curvature of $\gro_1$ and $c$ is a constant.  Since $\gro_1$ is extremal so is $\gro$, so the vector f\/ield $\partial^\#_{\gro} s_\gro$ is holomorphic. Moreover, this lifts to a transversely holomorphic vector f\/ield on $Y^{p,q}$ which equals $\partial^\#_gs_g$ where $s_g$ is the scalar curvature of the Sasakian metric~$g$ implying that $g$ is extremal.
\end{proof}

Again applying the openness theorem \cite{BGS06}, we obtain an open cone of extremal metrics containing the extremal structure of Proposition~\ref{otherext}. At this time we do not know whether the full extremal set $\ge_3(q)$ in $\gt_3^+(q)$ is connected. Note also that the pair $(\bp_-,\bp_+)$ is never $(1,1)$, so the f\/irst factor is always a weighted projective line, that is it is $\bbc\bbp^1$ with a non-trivial orbifold structure.

In \cite{Boy10a,Boy10b} the author introduced the notion of a Sasaki bouquet $\gB_N(\cald)$ associated to a~contact structure which consists of a union of $N$ Sasaki cones. Thus, associated to the contact structure~$\cald_p$ on $S^2\times S^3$ we have a bouquet $\gB_{\phi(p)}(\cald_p)$ consisting of $\phi(p)$ Sasaki cones. Moreover, each Sasaki cone has an extremal subset. Exactly how large these extremal sets are is unknown at this time. The only toric contact structures where the Sasakian extremal set is known to f\/ill up the Sasaki cone is for the standard toric contact structure on the sphere~$S^{2n+1}$~\cite{BGS06} and the standard completely integrable contact structure of Reeb type $(\gH^{2n+1},\cald,\eta,1)$ of Example~\ref{Heis} on the Heisenberg group~\cite{Boy09}.

\subsection*{Acknowledgements}

{\sloppy During the conference I enjoyed conversations with E.~Kalnins, N.~Kamran, J.~Kress, W.~Miller~Jr., and P.~Winternitz. I also want to thank J.~Pati, my collaborator in~\cite{BoPa10} without whom the present paper could not have been written.

}

\pdfbookmark[1]{References}{ref}
\LastPageEnding

\end{document}